\documentclass[reqno]{amsart}
\usepackage{amssymb}
\usepackage{hyperref}
\usepackage[sort,nocompress]{cite} 
\newtheorem{theorem}{Theorem}[section]
\newtheorem{lemma}[theorem]{Lemma}
\usepackage{amsaddr}
\usepackage{graphicx}
\usepackage{geometry}
\theoremstyle{definition}

\usepackage{xcolor}
\theoremstyle{remark}
\newtheorem{remark}[theorem]{Remark}

\numberwithin{equation}{section}
\makeatletter
\renewcommand\subsection{\@startsection{subsection}{2}%
	\z@{-.5\linespacing\@plus-.7\linespacing}{.5\linespacing}%
	{\normalfont\scshape}}
\renewcommand\subsubsection{\@startsection{subsubsection}{3}%
	\z@{.5\linespacing\@plus.7\linespacing}{-.5em}%
	{\normalfont\scshape}}
\makeatother

\begin{document}
	
	\title[{ Numerical Acoustic Wave Equation}]{Efficient and stable finite difference modelling of acoustic wave propagation in variable-density media  }
	
	\author{Da Li, Keran Li, Wenyuan Liao*}
	\address{Department of Mathematics and Statistics\\
		University of Calgary, AB, T2N 1N4, Canada}
	\email{da.li1@ucalgary.ca}
	
	\email{keran.li1@ucalgary.ca}
	
	\email{wliao@ucalgary.ca}
	
	
	
	
	
	\begin{abstract}
		In this paper, we consider the development and analysis of a new explicit compact high-order finite difference scheme for acoustic wave equation formulated in divergence form, which is widely used to describe seismic wave propagation through a heterogeneous media with variable media density and acoustic velocity. The new scheme is compact and of  fourth-order accuracy in space and second-order accuracy in time. The compactness of the scheme is obtained by the so-called combined finite difference method, which utilizes the boundary values of the spatial derivatives and those boundary values are obtained by one-sided finite difference approximation. An empirical stability analysis has been conducted to obtain the Currant-Friedrichs-Lewy (CFL) condition, which confirmed the conditional stability of the new scheme. Four numerical examples have been conducted  to validate the convergence and effectiveness of the new scheme. The application of the new scheme to a realistic wave propagation problem with Perfect Matched Layer boundary condition is also validated in this paper as well.
		
		\smallskip
		\noindent \textbf{Keywords:} Acoustic Wave Equation, Heterogeneous Media, High-order Compact Finite Difference Scheme, Stability.
	\end{abstract}
	
	\maketitle

	\section{Introduction}
	
	Finite difference method is one of the most popular numerical methods to solve partial differential equations for some reasons, for example, easy implementation and high efficiency. In particular, high-order (high-order here means that the spatial accuracy order is greater than or equal to four) finite difference methods have attracted increasing attentions recently from researchers in both science and engineering. Due to the great deal of efforts by many researchers and the increasing demand for highly accurate numerical solutions from the geophysics community, many finite difference methods have been developed to solve acoustic wave equations \cite{alford1974accuracy,liu2009implicit,tam1993dispersion,yang2012central,Mattsson2006,Bilbao2018}.  In addition to the high accuracy, there are many  other benefits from using high-order numerical method, such as the effectiveness in minimizing dispersion errors \cite{finkelstein2007finite},  fewer grid points being required than the conventional finite difference methods \cite{etgen2007computational}, etc.  Furthermore, it has been pointed out that high-order finite difference methods allow a more coarse sampling rate \cite{levander1988fourth}. 
	
	However, many high-order finite difference methods are not compact, which usually causes difficulties in treating the boundary conditions. For example, in 1D case, the conventional 3-point stencil can be used to approximate the first derivatives or second derivatives with at most 2nd-order accuracy. For higher-order accuracy, more points will be involved, thus,  it will be difficult to implement the method if  only one layer of boundary condition has been specified. Compact high-order finite difference methods are introduced to address those difficulties. In \cite{lele1992compact}, the author proposed a family of compact finite difference approximations to the first and second derivatives which resulted in 4th-order approximations to the first and second derivatives involving three points only in 1D case. This approach, although being very efficient and compact, needs boundary conditions for $u_x$ and $u_{xx}$, which are usually not available in the original problem. Under the smoothness assumption, a wave-equation based approach was proposed \cite{li2020efficient} to approximate the boundary conditions  of $u_{xx}$ with arbitrary  high order accuracy.  It is worth to mention that many high-order compact finite difference methods are only available for constant coefficients cases and fail in variable coefficients cases. To address this issue, new schemes have been proposed. Recently, two 4th-order compact finite difference schemes for acoustic wave equations  with variable sound speed was introduced in \cite{li2020efficient,li2019compact}. In \cite{britt2018high}, an energy norm based method was applied to analyze the stability of the finite difference methods for acoustic wave equation with variable coefficients. However the self-adjointness of the discretized Laplacian operator is required. Some of other works on compact high-order finite difference methods can be found in \cite{kim1996optimized,liao2014dispersion}.   The authors of \cite{shukla2005derivation} also discussed  some of the recent compact finite difference methods.
	
	The methods proposed in \cite{li2019compact,li2020efficient} work very effectively for variable acoustic velocity with constant media density cases, however, many real-world applications deal with the variable media density case, for which the Laplacian operator in spatial dimension is replaced by
	the divergence form given by $\nabla \cdot \left(\frac{1}{\rho} \nabla u\right)$.  It has also been shown \cite{cohen1996} that  the fourth-order accuracy can be obtained only if the density is smooth. 
	
	To overcome the forementioned drawback, a novel compact finite difference scheme for acoustic wave equation in heterogeneous media is proposed and investigated in this paper. The new scheme is of 4th-order in space and 2nd-order in time. Due to the heterogeneity of the media density, the regular Laplacian that appears in the wave equation  is replaced by  Laplacian in divergence form $\nabla \cdot \left(\frac{1}{\rho}\nabla\right)$, where $\rho$ indicates the density  of the media.  This scheme needs the values of the derivatives of the solution on the boundary, which can be approximated by Dirichlet boundary conditions and one-sided finite difference approximation.  The rest part of this paper is organized as the follows. In Section 2, the new scheme is proposed. An empirical stability condition is obtained by analyzing the eigenvalues of the discretized differential operators in Section 3. Numerical examples are solved in Section 4, which also includes the application of the new scheme to acoustic wave equation in heterogeneous media with Perfect Matched Layer (PML). 
	Finally, conclusions and some possible future works are addressed in Section 5.
	The derivation of the acoustic wave equation with PML boundary conditions is also included in Appendix A.
	
	\section{The new compact scheme}

	In this paper, the following 3D acoustic wave equation in heterogeneous media with initial condition and Dirichlet boundary conditions is considered,
	
	\begin{equation}\label{first eq.}
	\frac{1}{\rho c^2}u_{tt} - \nabla \cdot \left(\frac{1}{\rho} \nabla u\right) = s,\ (t,x,y,z)\in[0,T]\times\Omega,
	\end{equation} 
	where $s$ is the source term, $c$ the acoustic velocity and $\rho$ the media density. Here $s$ is a function of $x,y,z$ and t, while  $c$ and $\rho$ are functions of $x$, $y$ and $z$. The solution $u$ represents the wave pressure.
	
	Note that 
	\begin{equation}
	\nabla \cdot \left(\frac{1}{\rho} \nabla u\right) = \partial_x \left(\frac{1}{\rho } \partial_x u\right) + \partial_y \left(\frac{1}{\rho } \partial_y u\right) + \partial_z \left(\frac{1}{\rho } \partial_z u\right),
	\end{equation}
	it is necessary to approximate each term of the right-hand side with 4th-order accuracy to obtain an overall 4th-order spatial accuracy.
	
	Consider a single-variable function, $v(x)$. In \cite{lele1992compact, chu1998three}, the authors proposed the so-called combined compact 4th-order finite difference approximation for the first derivative, 
	\begin{equation}\label{1st derivative}
	\frac{1}{4}v'_{i-1} + v'_i + \frac{1}{4}v'_{i+1} = \frac{1}{h}\left(-\frac{3}{4}v_{i-1} + \frac{3}{4}v_{i+1} \right)
	\end{equation}
	where $v_{i}$ is the function $v$ evaluated at the grid point $x_i$, which is $v(x_i)$, $h$ the grid size. Similarly, the second derivative can be approximated in 4th-order by
	\begin{equation}\label{2nd derivative}
	\frac{1}{10}v''_{i-1} + v''_i + \frac{1}{10}v''_{i+1} = \frac{1}{h^2}\left(\frac{6}{5}v_{i-1} - \frac{12}{5}v_i + \frac{6}{5}v_{i+1} \right).
	\end{equation}
	The 4th-order accuracy of the scheme (\ref{1st derivative})(\ref{2nd derivative}) can be verified by Taylor expansion.
	
	
	For simplicity, we assume that $\Omega$ is a rectangle box defined as
	\begin{equation*}
		\Omega = [x_{min},x_{max}]\times[y_{min},y_{max}]\times[z_{min},z_{max}],
	\end{equation*}
	which is discretized into an $(N_x+2) \times (N_y+2) \times (N_z+2)$ grid with spatial grid sizes $h_x = \dfrac{x_{max} - x_{min}}{N_x + 1}$, $h_y = \dfrac{y_{max} - y_{min}}{N_y +1}$ and $h_z = \dfrac{z_{max} - z_{min}}{N_z + 1}$. Then the initial-boundary value problem of the 3D acoustic wave equation can be rewritten in this form
	\begin{equation}\label{acoustic eq}
	\begin{cases}
	\frac{1}{\rho(x,y,z) c^2(x,y,z)}u_{tt} - \nabla \cdot \left(\frac{1}{\rho(x,y,z)} \nabla u\right) = s(t,x,y,z),\\
	u|_{t=0} = \alpha(x,y,z),\ u_t|_{t=0} = \beta(x,y,z),\\
	u|_{x=x_{min}} = f_0(t,y,z),\ u|_{x=x_{max}} = f_1(t,y,z), \\
	u|_{y=y_{min}} = g_0(t,x,z),\ u|_{y=y_{max}} = g_1(t,x,z), \\
	u|_{z=z_{min}} = h_0(t,x,y),\ u|_{z=z_{max}} = h_1(t,x,y). \\
	\end{cases}
	\end{equation}
	Denoted by $\tau$ the time step, $u_{i,j,k}^n$ the numerical solution at grid point $(x_i,y_j,z_k)=(x_{min} +i h_x,y_{min}+j h_y,z_{min}+ k h_z)$ and time level $t_n = n \tau$, also define $\rho_{i,j,k}$, $c_{i,j,k}$ and $s^n_{i,j,k}$ in similar way. Define the following vectors
	\begin{equation}
	u_{*,j,k}^n =
	\begin{pmatrix}
	u_{1,j,k}^n\\
	u_{2,j,k}^n\\
	\vdots\\
	u_{N_x,j,k}^n
	\end{pmatrix}_{N_x \times 1},\ 
	(u_x)_{*,j,k}^n =
	\begin{pmatrix}
	(u_x)_{1,j,k}^n\\
	(u_x)_{2,j,k}^n\\
	\vdots\\
	(u_x)_{N_x,j,k}^n
	\end{pmatrix}_{N_x \times 1}
	\end{equation}
	and
	\begin{equation}
	\left[\left(\frac{1}{\rho}u_x\right)_x\right]_{*,j,k}^n =
	\begin{pmatrix}
	\left[\left(\frac{1}{\rho}u_x\right)_x\right]_{1,j,k}^n\\
	\left[\left(\frac{1}{\rho}u_x\right)_x\right]_{2,j,k}^n\\
	\vdots\\
	\left[\left(\frac{1}{\rho}u_x\right)_x\right]_{N_x,j,k}^n
	\end{pmatrix}_{N_x \times 1},
	\end{equation}
	where $1\leqslant i\leqslant N_x$, $1\leqslant j\leqslant N_y$, $1\leqslant k\leqslant N_z$. Throughout the rest of this paper, a quantity represented by a  lower-case letter with $*$ being one of its subscripts  denotes the vector form of the quantity, which is defined in the similar way as above. Then from (\ref{1st derivative}) one has the approximation equations for $\left[\left(\frac{1}{\rho}u_x\right)_x\right]_{*,j,k}^n$
	\begin{equation}\label{u_x}
	A_x (u_{x})_{*,j,k}^n+\frac{1}{4}w_x^{a,n} = \frac{1}{h_x} \left(B_x u_{*,j,k}^n+\frac{3}{4}w_x^{b,n} \right),
	\end{equation}
	and
	\begin{equation}\label{u_xx}
	A_x \left[\left(\frac{1}{\rho}u_x\right)_x\right]_{*,j,k}^n+\frac{1}{4}w_{xx}^{a,n} = \frac{1}{h_x} \left(B_x \left[\frac{1}{\rho}u_x\right]_{*,j,k}^n+\frac{3}{4}w_{xx}^{b,n} \right),
	\end{equation}
	where
	\begin{equation}
	A_x = 
	\begin{pmatrix}
	1 & \frac{1}{4} \\
	\frac{1}{4} & 1 &\frac{1}{4} \\
	& \dots &\dots\\
	& \frac{1}{4} & 1 &\frac{1}{4} \\
	& & \frac{1}{4} & 1 \\
	\end{pmatrix}_{N_x \times N_x},
	B_x = 
	\begin{pmatrix}
	0 & \frac{3}{4} \\
	-\frac{3}{4} & 0 & \frac{3}{4} \\
	& \cdots & \cdots \\
	& -\frac{3}{4} & 0 & \frac{3}{4} \\
	& & -\frac{3}{4} & 0
	\end{pmatrix}_{N_x \times N_x}
	\end{equation}
	are tridiagonal matrices, and
	\begin{equation}
	w_x^{a,n} =
	\begin{pmatrix}
	(u_{x})_{0,j,k}^n \\
	0\\
	\vdots\\
	0\\
	(u_{x})_{N_x+1,j,k}^n
	\end{pmatrix}_{N_x \times 1},\ 
	w_x^{b,n} =
	\begin{pmatrix}
	-u_{0,j,k}^n \\
	0\\
	\vdots\\
	0\\
	u_{N_x+1,j,k}^n
	\end{pmatrix}_{N_x \times 1}
	\end{equation}
	\begin{equation}
	w_{xx}^{a,n} =
	\begin{pmatrix}
	\left[\left(\frac{1}{\rho}u_x\right)_x\right]_{0,j,k}^n \\
	0\\
	\vdots\\
	0\\
	\left[\left(\frac{1}{\rho}u_x\right)_x\right]_{N_x+1,j,k}^n
	\end{pmatrix}_{N_x \times 1},\ 
	w_{xx}^{b,n} =
	\begin{pmatrix}
	-\left[\frac{1}{\rho}u_x\right]_{0,j,k}^n \\
	0\\
	\vdots\\
	0\\
	\left[\frac{1}{\rho}u_x\right]_{N_x+1,j,k}^n
	\end{pmatrix}_{N_x \times 1}
	\end{equation}
	represent the boundary values. Note that the negative sign of the first components of $w_x^{b,n}$ and $w_{xx}^{b,n}$ results from the negative term of the right-hand side of (\ref{1st derivative}). The approximation equations for $\left[\left(\frac{1}{\rho}u_y\right)_y\right]_{i,*,k}^n$ 
	and $\left[\left(\frac{1}{\rho}u_z\right)_z\right]_{i,j,*}^n$ can be obtained similarly. Note that for the boundary values vectors $w_x$'s and $w_{xx}$'s, only $w_x^{b,n}$ can be evaluated directly from the boundary conditions of the equation. For the boundary values occur in $w_x^{a,n}$, $w_{xx}^{b,n}$ and $w_{xx}^{a,n}$, they can be obtained by the 4th-order one-sided finite difference approximation for the first derivatives,
	\begin{equation}\label{sided start}
	v'_0 = \frac{1}{h_x}\left(
	-\frac{25}{12}v_0 + 4 v_1 -3 v_2 +\frac{4}{3} v_3 - \frac{1}{4}v_4
	\right),
	\end{equation}
	and
	\begin{equation}\label{sided end}
	v'_{N_x+1} = \frac{1}{h_x}\left(
	\frac{25}{12}v_{N_x+1} - 4 v_{N_x} +3 v_{N_x-1} -\frac{4}{3} v_{N_x-2} + \frac{1}{4}v_{N_x-3}
	\right).
	\end{equation}
	Suppose that  $u^n_{*,j,k}$ is known, the following steps show how to obtain $\left[\left(\frac{1}{\rho}u_x\right)_x\right]_{*,j,k}^n$ from the above discussion.
	\begin{enumerate}
		\item Both $u^n_{0,j,k}$ and $u^n_{N_x+1,j,k}$ are evaluated from the boundary conditions of the equation, thus $w_x^{b,n}$ can be obtained.
		\item Use (\ref{sided start})(\ref{sided end}) to approximate $(u_x)^n_{0,j,k}$ and $(u_x)^n_{N_x+1,j,k}$ by $u^n_{*,j,k}$, thus $w_x^{a,n}$ will be known.
		\item Solve (\ref{u_x}) to obtain $(u_x)^n_{*,j,k}$, then $\left(\frac{1}{\rho}u_x\right)^n_{*,j,k}$ will be known.
		\item Use (\ref{sided start})(\ref{sided end}) to approximate $\left(\frac{1}{\rho}u_x\right)^n_{0,j,k}$ and $\left(\frac{1}{\rho}u_x\right)^n_{N_x+1,j,k}$ by $\left(\frac{1}{\rho}u_x\right)^n_{*,j,k}$, thus $w_{xx}^{b,n}$ will be known.
		\item Use (\ref{sided start})(\ref{sided end}) to approximate $\left[\left(\frac{1}{\rho}u_x\right)_x\right]_{0,j,k}^n$ and $\left[\left(\frac{1}{\rho}u_x\right)_x\right]_{N_x+1,j,k}^n$ by $\left(\frac{1}{\rho}u_x\right)^n_{0,j,k}$, $\left(\frac{1}{\rho}u_x\right)^n_{N_x+1,j,k}$ and $\left(\frac{1}{\rho}u_x\right)^n_{*,j,k}$, thus $w_{xx}^{a,n}$ will be known.
		\item Solve (\ref{u_xx}) to obtain $\left[(\frac{1}{\rho}u_x)_x\right]^n_{*,j,k}$.
	\end{enumerate}
	The assumption that $u^n_{*,j,k}$ is known is reasonable, since this paper uses the conventional 2nd-order central finite difference to approximate the second derivative $u_{tt}$, i.e. $u^{n+1}$ is solved from $u^{n}$ and $u^{n-1}$. The derivative terms $\left[\left(\frac{1}{\rho}u_y\right)_y\right]^n_{i,*,k}$ and $\left[\left(\frac{1}{\rho}u_z\right)_z\right]^n_{i,j,*}$ can be obtained similarly.
	
	\begin{remark}
		If the media density $\rho$ is differentiable, then one can also consider the equivalent form of the Laplacian $\nabla \cdot (\frac{1}{\rho}\nabla u) = \nabla \frac{1}{\rho} \cdot \nabla u + \frac{1}{\rho}\Delta u$, which will be simpler to implement, where $\nabla u$ is approximated by the 4th-order compact finite difference scheme as above, and $\Delta u$ can be approximated by the 4th-order compact finite difference in \cite{li2020efficient}. Note that the approximation of the second spatial derivatives of $u$ on the boundary requires a one-sided finite difference approximation similar to (\ref{sided start}) and (\ref{sided end}) with different coefficients.
	\end{remark}
	
	To initialize the solving process, one also needs $u^{-1}_{i,j,k} = u(-\tau,x_i,y_j,z_k)$, which can be obtained by
	\begin{equation}\label{u at -tau}
	\begin{split}
	u^{-1}_{i,j,k} =& u^0_{i,j,k} - \tau (u_t)^0_{i,j,k} + \frac{1}{2}{\tau}^2 (u_{tt})^0_{i,j,k} - \frac{1}{6}{\tau}^3 (u_{ttt})^0_{i,j,k} + O(\tau^4)\\
	=& \alpha_{i,j,k} - \tau \beta_{i,j,k} + \frac{1}{2}{\tau}^2\left\{\left(\rho c^2\right)_{i,j,k}\left[\nabla \cdot \left(\frac{1}{\rho} \nabla \alpha\right) \right]_{i,j,k}+s^0_{i,j,k}\right\} \\
	-&\frac{1}{6}{\tau}^3\left\{\left(\rho c^2\right)_{i,j,k}\left[\nabla \cdot \left(\frac{1}{\rho} \nabla \beta\right) \right]_{i,j,k}+(\partial_t s)^0_{i,j,k}\right\}+O(\tau^4).
	\end{split},
	\end{equation}
	where $\alpha = u|_{t = 0}$ and $\beta = u_t|_{t = 0}$ are the initial conditions.
	
	Finally, a compact finite difference scheme with error $O(\tau^2)+O(h_x^4)+O(h_y^4)+O(h_z^4)$ is obtained as
	\begin{equation}
	u^{n+1}_{i,j,k} = \tau^2 \left\{\left(\rho c^2\right)_{i,j,k}\left[\nabla \cdot \left(\frac{1}{\rho} \nabla u\right) \right]^n_{i,j,k}+s^n_{i,j,k} \right\}+2 u^{n}_{i,j,k} - u^{n-1}_{i,j,k}, \ n = 0, 1,2,\cdots
	\end{equation}
	with $\left[\nabla \cdot \left(\frac{1}{\rho} \nabla u\right) \right]^n_{i,j,k}$ obtained from $\left[(\frac{1}{\rho}u_x)_x\right]^n_{*,j,k}$,  $\left[(\frac{1}{\rho}u_y)_y\right]^n_{i,*,k}$ and $\left[(\frac{1}{\rho}u_z)_z\right]^n_{i,j,*}$, and $u^{-1}_{i,j,k}$ from (\ref{u at -tau}).
	
	\section{Stability Analysis of the New Scheme}
	
	Consider the acoustic wave equation with zero boundary conditions and zero source term
	\begin{equation}
	\frac{1}{\rho c^2}u_{tt} - \nabla \cdot \left(\frac{1}{\rho} \nabla u\right) = 0,
	\end{equation}
	or
	\begin{equation}\label{ZeroSourceBoundary}
	u_{tt} -  c^2 \rho \left[\nabla \cdot \left(\frac{1}{\rho} \nabla u\right)\right] = 0.
	\end{equation}
	For simplicity assume $h_x = h_y = h_z = h$, $N_x = N_y = N_z = N$, $\tau$ the time step size. Also let $U^n$ be the vector form of the numerical solution $u^n_{i,j,k}$ denoted by
	\begin{equation}
	U^n =
	\begin{pmatrix}
	u^n_{1,1,1} & \dots & u^n_{N,1,1} &u^n_{1,2,1} & \dots & u^n_{N,2,1} & u^n_{1,3,1} & \dots & u^n_{N,3,1} & \dots
	\end{pmatrix}^T
	\end{equation}
	i.e. $U^n$ is an $N^3\times 1$ vector in which $u^n_{i,j,k}$ is located at the $(k N^2 + j N + i)$-th row. Also define $D_x^2 U$, $D_y^2 U$ and $D_z^2 U$ as the vector form of the derivative term $\left[\left(\frac{1}{\rho}u_x\right)_x\right]^n_{i,j,k}$, $\left[\left(\frac{1}{\rho}u_y\right)_y\right]^n_{i,j,k}$ and $\left[\left(\frac{1}{\rho}u_z\right)_z\right]^n_{i,j,k}$, respectively. 
	Then equation (\ref{ZeroSourceBoundary}) can be approximated by
	\begin{equation}\label{ZeroSourceBoundary D}
	\frac{1}{\tau^2}\delta_t^2 U^n - C^2 Q\left(D_x^2 U^n + D_y^2 U^n + D_z^2 U^n\right) = 0,
	\end{equation}
	where $Q$ is the diagonal matrix for the function $\rho$
	\begin{equation}
	Q = diag
	\begin{pmatrix}
	q_{1,1,1} & \dots & q_{N,1,1} &q_{1,2,1} & \dots & q_{N,2,1} & q_{1,3,1} & \dots & q_{N,3,1} & \dots
	\end{pmatrix}
	\end{equation}
	and $C$ is the diagonal matrix for the function $c$
	\begin{equation}
	C = diag
	\begin{pmatrix}
	c_{1,1,1} & \dots & c_{N,1,1} &c_{1,2,1} & \dots & c_{N,2,1} & c_{1,3,1} & \dots & c_{N,3,1} & \dots
	\end{pmatrix}.
	\end{equation}
	In other words, $Q$ and $C$ are $N^3 \times N^3$ diagonal matrices whose diagonal entries are the vector form of $q_{i,j,k}$ and $c_{i,j,k}$, respectively.
	
	Let
	\begin{equation}
	A = 
	\begin{pmatrix}
	1 & \frac{1}{4} \\
	\frac{1}{4} & 1 &\frac{1}{4} \\
	& \dots &\dots\\
	& \frac{1}{4} & 1 &\frac{1}{4} \\
	& & \frac{1}{4} & 1 \\
	\end{pmatrix}_{N\times N},
	B = 
	\begin{pmatrix}
	0 & \frac{3}{4} \\
	-\frac{3}{4} & 0 & \frac{3}{4} \\
	& \cdots & \cdots \\
	& -\frac{3}{4} & 0 & \frac{3}{4} \\
	& & -\frac{3}{4} & 0
	\end{pmatrix}_{N\times N}
	\end{equation}
	Then the spatial derivative terms $D_x^2 U$, $D_y^2 U$ and $D_z^2 U$ of equation (\ref{ZeroSourceBoundary D}) can be written as
	\begin{equation}
	D_x^2 U^n = \frac{1}{h^2}\left[A_1^{-1}B_1Q^{-1}A_1^{-1}B_1\right] U^n,
	\end{equation}
	\begin{equation}
	D_y^2 U^n = \frac{1}{h^2}\left[A_2^{-1}B_2Q^{-1}A_2^{-1}B_2\right] U^n,
	\end{equation}
	\begin{equation}
	D_z^2 U^n = \frac{1}{h^2}\left[A_3^{-1}B_3Q^{-1}A_3^{-1}B_3\right] U^n,
	\end{equation}
	where
	\begin{equation}
	A_1 = A\otimes I_N \otimes I_N,\ B_1 = B\otimes I_N \otimes I_N,
	\end{equation}
	\begin{equation}
	A_2 = I_N\otimes A \otimes I_N,\ B_2 = I_N\otimes B \otimes I_N,
	\end{equation}
	\begin{equation}
	A_3 = I_N\otimes I_N \otimes A,\ B_3 = I_N\otimes I_N \otimes B,
	\end{equation}
	with $\otimes$ being  the Kronecker product and $I_N$ the $N\times N$ identity matrix.
	
	Now define
	\begin{equation}
	L = C^2 Q \left(A_1^{-1}B_1Q^{-1}A_1^{-1}B_1 + A_2^{-1}B_2Q^{-1}A_2^{-1}B_2 + A_3^{-1}B_3Q^{-1}A_3^{-1}B_3\right),
	\end{equation}
	then the discretized equation (\ref{ZeroSourceBoundary D}) can be written as
	\begin{equation}\label{L eq original}
	\frac{1}{\tau^2}\delta_t^2 U^n - \frac{1}{h^2} L U^n = 0,
	\end{equation}
	where $\delta_t^2$ is the conventional temporal 2nd-order central finite difference operator
	\begin{equation}
	\delta_t^2 U^n = U^{n+1} - 2U^{n} + U^{n-1}.
	\end{equation}
	
	It will be very useful if the estimate of the eigenvalues of the spatial difference operator is known. However, it is very difficult to obtain the estimate for $L$,  due to the variant coefficients in $Q$ and $C$. It is empirical that the eigenvalue of $L$ with largest absolute value determines the stability of a scheme for wave equations. Thus,  consider freezing the coefficients of $L$, which  leads to a variant of $L$,
	\begin{equation}\label{L tilde}
	\tilde{L} = c_{max}^2 \frac{q_{max}}{q_{min}} \left(A_1^{-1}B_1 A_1^{-1}B_1 + A_2^{-1}B_2 A_2^{-1}B_2 + A_3^{-1}B_3 A_3^{-1}B_3\right),
	\end{equation}
	where $c_{max}$ and $q_{max}$ are the maximum entries of $C$ and $Q$, respectively, $q_{min}$ the minimum entry of $Q$. Note that all of the entries of $C$ and $Q$ are positive. 
	
	\subsection{Estimation of the Eigenvalues}
	
	To estimate the spectrum of $\tilde{L}$, we introduce the following lemmas, which can be found in \cite{horn1991topics}.
	
	\begin{lemma}\label{Kronecker property}
		Kronecker product is associative. For square matrices $K$, $L$, $G$ and $H$, the following identities hold
		\begin{equation}
		\begin{split}
		I_M \otimes I_N &= I_{MN},\\
		(K\otimes J)(G\otimes H) &= (KJ)\otimes (GH),\\
		(K\otimes J)^{-1}&= (K^{-1}\otimes J^{-1}),\\
		(K\otimes J)^{T}&= (K^{T}\otimes J^{T}).
		\end{split}
		\end{equation}
	\end{lemma}
	\begin{lemma}\label{Kronecker eig}
		For square matrices $K$ and $J$, any eigenvalue of $K \otimes J$ arises as a product of eigenvalues of $K$ and $J$. If $\lambda_K$ and $\lambda_J$ are eigenvalues of $K$ and $J$, respectively, then $\lambda_K \cdot \lambda_J$ is an eigenvalue  of $K \otimes J$.
	\end{lemma}
	With the above lemmas, one has
	\begin{equation}\label{ABAB}
	\begin{split}
	A_1^{-1}B_1 A_1^{-1}B_1 = (A^{-1}B A^{-1}B)\otimes I_N \otimes I_N, \\
	A_2^{-1}B_2 A_2^{-1}B_2 =I_N \otimes(A^{-1}B A^{-1}B)\otimes  I_N, \\
	A_3^{-1}B_3 A_3^{-1}B_3 =I_N\otimes I_N \otimes  (A^{-1}B A^{-1}B),
	\end{split}
	\end{equation}
	thus
	\begin{equation}
	\sigma\left(A_1^{-1}B_1 A_1^{-1}B_1 \right) = \sigma \left(A_2^{-1}B_2 A_2^{-1}B_2\right) = \sigma \left(A_3^{-1}B_3 A_3^{-1}B_3 \right) = \sigma \left( A^{-1}B A^{-1}B\right).
	\end{equation}
	Recall the definition of Kronecker sum of two square matrices $K_{M\times M}$ and $L_{N\times N}$
	\begin{equation}
	K \oplus J = (I_N \otimes K) + (J \otimes I_M).
	\end{equation}
	The following lemma can also be found in \cite{horn1991topics}.
	\begin{lemma}
		For square matrices $K$ and $J$, any eigenvalue of $K \oplus J$ arises as a sum of eigenvalues of $K$ and $J$. If $\lambda_K$ and $\lambda_J$ are eigenvalues of $K$ and $J$, respectively, then $\lambda_K+\lambda_J$ is an eigenvalue  of $K\oplus J$.	
	\end{lemma}
	Base on that, we can verify that
	\begin{equation}\label{KJG}
	K \oplus J \oplus G = (I_N \otimes I_N \otimes K) + (I_N \otimes J \otimes I_N) + (G \otimes I_N \otimes I_N)
	\end{equation}
	and  
	\begin{equation}\label{KJG e.v.}
	\sigma(K \oplus J \oplus G) = \{\lambda_K + \lambda_J +\lambda_G,\ \lambda_K \in \sigma(K),\ \lambda_J \in \sigma(J),\ \lambda_G \in \sigma(G)
	\}.
	\end{equation}
	With (\ref{L tilde})(\ref{ABAB}), one has
	\begin{equation}\label{L tilde sum}
	\tilde{L} = c_{max}^2 \frac{q_{max}}{q_{min}}\left[(A^{-1}B A^{-1}B) \oplus (A^{-1}B A^{-1}B) \oplus (A^{-1}B A^{-1}B) \right].
	\end{equation}
	Now it is necessary to find out the eigenvalues of $A^{-1}B A^{-1}B$. Since $A$ and $B$ are tridiagonal Toeplitz matrices, their eigenvalues are given by
	\begin{equation}
	\sigma(A) = \left\{ 1+\frac{1}{2}\cos \frac{l\pi}{N+1},\ 1\leqslant l \leqslant N \right\} \subset \left(\frac{1}{2},\frac{3}{2}\right),
	\end{equation} 
	therefore, $A$ is a positive definite matrix and
	\begin{equation}
	\sigma(B) =\left \{\frac{3\sqrt{-1}}{2}\cos \frac{l\pi}{N+1},\ 1\leqslant l \leqslant N\right \}\subset \left(-\frac{3\sqrt{-1}}{2},\frac{3\sqrt{-1}}{2}\right).
	\end{equation} 
	
	Recall that for a square matrix $K$, the numerical range of $F(K)$ is defined as
	\begin{equation}
	F(K) = \left\{ \frac{\vec{v}^T K \vec{v}}{\vec{v}^T \vec{v}}, v\ \mbox{non-zero complex vectors} \right\},
	\end{equation}
	which is a closed bounded convex set containing all of the eigenvalues of $K$. Furthermore, if $K$ is a normal matrix, $F(K)$ is exactly the convex closure of all the eigenvalues of $K$. 
	\begin{theorem}
		For two real square matrices $K$ and $J$ of the same size, if $K$ is a symmetric semi-positive definite matrix and $\lambda$ is an eigenvalue of $KJ$, then 
		\begin{equation}
		\lambda \in F(K)F(J) = \{ \lambda_K\cdot \lambda_J,\ \lambda_K\in F(K)\ \mbox{and}\  \lambda_J \in F(J) \},
		\end{equation}
		in other words, $\sigma(KL)\subset F(K)F(L)$
	\end{theorem}
	The theorem above was proven in \cite{wielandt1973eigenvalues}.
	
	Now since $A$ is a symmetric semi-positive definite matrix, so is $A^{-1}$. Since $A^{-1}$ is symmetric and $B$ is anti-symmetric, thus, they are both normal matrices. Consequently, the numerical ranges of  $F(A^{-1})$ and $F(B)$ coincide with the closure of all the eigenvalues of $A^{-1}$ and $B$, respectively, i.e.
	\begin{equation}
	F(A^{-1})= \left(-2,-\frac{2}{3}\right),\ F(B)=  \left(-\frac{3\sqrt{-1}}{2},\frac{3\sqrt{-1}}{2}\right).
	\end{equation}
	Then the eigenvalues of $A^{-1}B$ can be estimated by
	\begin{equation}
	\sigma(A^{-1}B) \subset F(A^{-1})F(B) = \left[ 0, 3\sqrt{-1}\right).
	\end{equation} 
	\begin{remark}
		It is worth to notice that $0$ is an eigenvalue of $B$ if and only if $N$ is an odd number. This is also true for $A^{-1}B$.
	\end{remark}
	Then one has
	\begin{equation}
	\sigma(A^{-1}BA^{-1}B) \subset \left[ 0, -9\right),\ N\ \mbox{odd},
	\end{equation}
	and
	\begin{equation}
	\sigma(A^{-1}BA^{-1}B) \subset \left( 0, -9\right),\ N\ \mbox{even}.
	\end{equation}
	Finally, by (\ref{KJG e.v.})(\ref{L tilde sum}), one has
	\begin{equation}
	\sigma(\tilde{L}) \subset \left[ 0, -27c_{max}^2 \frac{q_{max}}{q_{min}}\right),\ N\ \mbox{odd},
	\end{equation}
	and
	\begin{equation}
	\sigma(\tilde{L}) \subset \left( 0, -27c_{max}^2 \frac{q_{max}}{q_{min}}\right),\ N\ \mbox{even}.
	\end{equation}
	
	\subsection{Stability}
	We now derive  an empirical CFL condition for the new scheme. Consider the coefficient-frozen version of the equation (\ref{L eq original})
	\begin{equation}
	\delta_t^2 U^n - \frac{\tau^2}{h^2} \tilde{L} U^n = 0
	\end{equation}
	i.e.
	\begin{equation}
	U^{n+1} - \left( 2 + \frac{\tau^2}{h^2} \tilde{L} \right) U^n + U^{n-1} = 0.
	\end{equation}
	Inspired by von Neumann analysis, one considers the following character equation
	\begin{equation}\label{character eq}
	\lambda^2 - \left(2 - 27r\right)\lambda + 1 = 0
	\end{equation}
	where $ r = c_{max}^2 \frac{q_{max}}{q_{min}}\frac{\tau^2}{h^2}$. Here the finite difference operator $\tilde{L}$ is replaced by $-27c_{max}^2 \frac{q_{max}}{q_{min}}$, which is the upper bound of the eigenvalue of $\tilde{L}$ with largest absolute value. In von Neumann analysis, the scheme is stable if the roots of the character equation (\ref{character eq}) are two complex numbers, i.e. the discriminant should be negative
	\begin{equation}
	\Delta = (2-27r)^2 - 4 = 27r(27r-4)<0,
	\end{equation}
	thus one has
	\begin{equation}
	0 < r = c_{max}^2 \frac{q_{max}}{q_{min}}\frac{\tau^2}{h^2} < \frac{4}{27},
	\end{equation}
	i.e.
	\begin{equation}\label{CFL}
	c_{max} \sqrt{\frac{q_{max}}{q_{min}}}\cdot \frac{\tau}{h} < \frac{2}{3\sqrt{3}}.
	\end{equation}
	The new scheme will be stable if the CFL condition (\ref{CFL}) is satisfied.
	
	\section{Numerical Experiments}
	
	In this section, four numerical examples are solved to demonstrate the accuracy and efficiency of the new scheme. The first example is solved to validate that the new scheme is of 2nd-order accuracy in time and 4th-order accuracy in space. The second example solves a zero initial and boundary conditions problem with Ricker wavelet source. The third example validates the effectiveness and the accuracy of the new scheme for acoustic wave equation with PML boundary conditions. The fourth example considers a more realistic problem, in which the Marmousi 2 model is used in the simulation of the seismic wave propagation.
	
	\subsection{Example 1}\label{e.g.:1}
	
	This example validates that the new scheme is of 2nd-order accuracy in time and 4th-order accuracy in space. Consider the acoustic wave equation defined on the domain $[0,1]\times[0,1]\times[0,1]$ and $t \in [0,1]$,
	\begin{equation}
	u_{tt} - \rho c^2\nabla \cdot \left(\frac{1}{\rho} \nabla u\right) = \rho c^2s,
	\end{equation} 
	where
	\begin{equation}
	\rho = e^{(-x-y-z)/3},
	\end{equation}
	\begin{equation}
	c^2 = 1+\frac{1}{2}xyz,
	\end{equation}
	and
	\begin{equation}
	\begin{split}
	\rho c^2 s = & -\sin(t)\cos(x+2y+3z) \\
	+ & \sin(t)\left(1+\frac{1}{2}xyz\right)\left[14\cos(x+2y+3z) + 2\sin(x+2y+3z)\right]
	\end{split}
	\end{equation}
	with initial and boundary conditions compatible to the analytic solution
	\begin{equation}
	u = \sin(t)\cos(x+2y+3z).
	\end{equation}
	The grid sizes are given by $h_x = h_y = h_z = h$ and the time step is $\tau = h^2$. Thus it is sufficient to show that in this numerical experiment, the accuracy order is $O(\tau^2) + O(h^4) = O(h^4)$. The max errors with different $h$ are listed in Table \ref{table:1}, which clearly validated the desired convergence order  calculated as
	\begin{equation}
	\text{Conv. Order} = \frac{\log\big[E(h_1)/E(h_2)\big]}{\log(h_1/h_2)}.
	\end{equation}
	
	\begin{table}
		\centering
		\caption{Numerical errors in max norm in Example \ref{e.g.:1} with $\tau = h^2$. Note that the CFL condition requires $\frac{\tau}{h} < 0.1156$.}
		\begin{tabular}{|c|c|c|c|c|c|}
			\hline
			h & 1/10 & 1/16 & 1/20 & 1/24 & 1/32\\
			\hline
			$\tau$ & 1/100 & 1/256 & 1/400 & 1/576 & 1/1024\\
			\hline
			$E$ & 7.6115e-05 & 9.5211e-06 & 3.8419e-06 	& 1.7292e-06 & 5.0288e-07 \\
			\hline
			Conv. Order	& - & 4.4228 & 4.0671 & 4.3786 & 4.2931	\\
			\hline
		\end{tabular}
		
		\label{table:1}
	\end{table}

	\subsection{Example 2}\label{e.g.:2}
	In this  example we solve the wave equation with Ricker wavelet source. Consider the equation on the region $\Omega = [0,2]\times[0,2]\times[0,2]$
	\begin{equation}
	\frac{1}{\rho}u_{tt} - \nabla\cdot\left(\frac{1}{\rho}\nabla u\right) = s,
	\end{equation}
	with zero initial and boundary conditions $u|_{t=0} = 0$, $u_{t}|_{t=0} = 0$ and $u|_{\partial \Omega} = 0$, where $\rho = 2z^2 + 1 $ and $s$ is the Ricker wavelet source given by
	\begin{equation}
	s(t,x,y,z) = \delta (x-x_0,y-y_0,z-z_0)[1-2\pi^2 f^2_p (t-d_r)^2]e^{-\pi^2 f_p^2 (t-d_r)^2},
	\end{equation}
	with dominant frequency $f_p = 10$ and temporal delay $d_r = \frac{1}{2f_p}$. The source is placed at the centre of the region with $x_0=y_0=z_0=1$. The equation is solved by the new scheme with $h_x = h_y = h_z = h = \frac{1}{40}$ and $\tau = \frac{1}{400}$. The stability condition which requires $\frac{\tau}{h}<\frac{2}{9\sqrt{3}} \approx 0.1283$ is satisfied. 
	
\begin{figure}[h]
		\includegraphics[width=10cm]{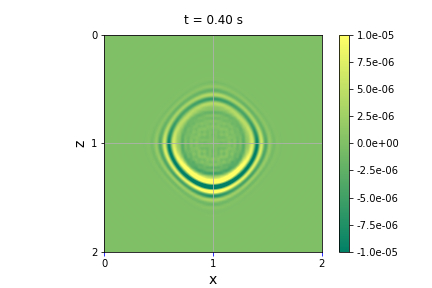}
		\caption{Snapshot of $y$-section at $t = 0.4s$ in Example \ref{e.g.:2}}.
		\label{fig:e.g. 2 0.4s}
	\end{figure}

\begin{figure}[h]
		\includegraphics[width=10cm]{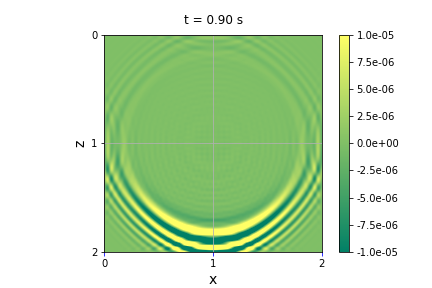}
		\caption{Snapshot of $y$-section at $t = 0.9s$ in Example \ref{e.g.:2}}.
		\label{fig:e.g. 2 0.9s}
	\end{figure}

\begin{figure}[h]
		\includegraphics[width=10cm]{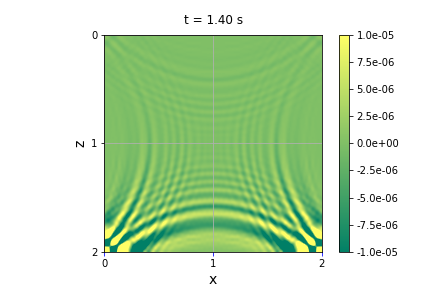}
		\caption{Snapshot of $y$-section at $t = 1.4s$ in Example. \ref{e.g.:2}}
		\label{fig:e.g. 2 1.4s}
	\end{figure}

	Three snapshots of $y$-section at $y=1$ at different times are plotted in Figure \ref{fig:e.g. 2 0.4s}, Figure \ref{fig:e.g. 2 0.9s} and Figure \ref{fig:e.g. 2 1.4s}. Note that the media density is given by $\rho = 2z^2 +1$, which depends only on $z$. Those figures show that the wave propagates at constant sound speed, as expected. However,  the energy is more concentrated in the area with higher media density. Finally, it is observed that the boundaries reflexes the wave back due to zero boundary conditions.
	\begin{remark}
		In real-world applications, variable media density usually results in variable acoustic velocity. However, in this example, for the sake of simplicity, the acoustic velocity is  normalized to a constant $c = 1$ to highlight the effect of the variable media density on the propagation of the acoustic wave.
	\end{remark}
	\begin{remark}
		It is worth to mention that Figure \ref{fig:e.g. 2 0.4s}, Figure \ref{fig:e.g. 2 0.9s} and Figure \ref{fig:e.g. 2 1.4s} are $y$-section of a wave in 3D space, thus the energy component of the wave in $y$-direction is not shown on those figures.
	\end{remark}
	
\begin{figure}[h]
		\includegraphics[width=10cm]{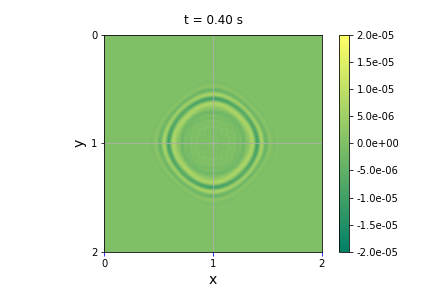}
		\caption{Snapshot of $z$-section at $t = 0.4s$ in Example \ref{e.g.:2}}.
		\label{fig:e.g. z2 0.4s}
	\end{figure}

\begin{figure}[h]
		\includegraphics[width=10cm]{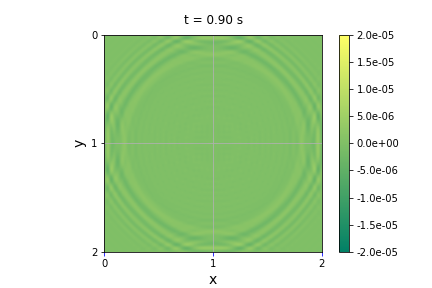}
		\caption{Snapshot of $z$-section at $t = 0.9s$ in Example \ref{e.g.:2}}.
		\label{fig:e.g. z2 0.9s}
	\end{figure}

\begin{figure}[h]
		\includegraphics[width=10cm]{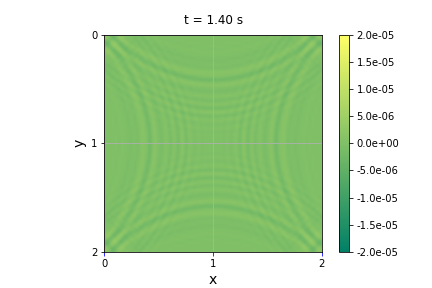}
		\caption{Snapshot of $z$-section at $t = 1.4s$ in Example. \ref{e.g.:2}}
		\label{fig:e.g. z2 1.4s}
	\end{figure}

	We then plot three snapshots of $z$-section at $z=1$ at different times in Figure \ref{fig:e.g. z2 0.4s}, Figure \ref{fig:e.g. z2 0.9s} and Figure \ref{fig:e.g. z2 1.4s}. It is clear that the wave propagates as if the density is constant, which is true as  the density is independent of $x$ and $y$. As can be seen, the energy is distributed symmetrically about the centre of the domain. Similarly, it is noticed that the boundaries reflexes the wave back due to zero boundary conditions.
	
	%
	%

	\subsection{Example 3}\label{e.g.:3}
	
	Perfectly Matched Layer (PML) is a technique to truncate the computational domain as simulating acoustic wave propagating in an unbounded domain, which is firstly introduced by Berenger \cite{berenger1994perfectly}.
	When PML is introduced,   the 3D acoustic wave equation (\ref{first eq.}) is modified to
	\begin{equation} \label{PML_eq}
	\begin{split} 
	\frac{1}{\rho c^2} \left(u_{tt} + \alpha u_t + \beta u + \gamma w \right) - \nabla\cdot \left(\frac{1}{\rho}\left(\vec v + \nabla u \right)\right) & = f, \\
	\vec v_t + H\vec v + J \nabla u - K \nabla w & = 0, \\
	w_t - u & = 0,
	\end{split}
	\end{equation}
	where $\sigma = \sigma_x + \sigma_y + \sigma_z$, $\zeta = \sigma_x\sigma_y + \sigma_x\sigma_z + \sigma_y\sigma_z$, $\gamma=\sigma_x\sigma_y\sigma_z$ and 
	\begin{equation*}
		H = \begin{pmatrix}
			\sigma_x & & \\ & \sigma_y & \\ & & \sigma_z
		\end{pmatrix},\ 
		K = \begin{pmatrix}
			\sigma_y\sigma_z & & \\ & \sigma_x\sigma_z & \\ & & \sigma_x\sigma_y
		\end{pmatrix}.
	\end{equation*}
	\begin{equation*}
		J = \begin{pmatrix}
			\sigma_x-\sigma_y-\sigma_z & & \\ & \sigma_y-\sigma_x-\sigma_z & \\ & & \sigma_z-\sigma_x-\sigma_y
		\end{pmatrix}.
	\end{equation*}
	The detailed derivation of the wave equation with PML is included in appendix \ref{appendix_A}.
	Here $\sigma_\lambda$ is a damping function with $\sigma_\iota = 0$ in the non-absorbing domain and $\sigma_\iota \neq 0$ in the absorbing layer for $\iota = x, y, z$.
	Also $\sigma_\iota$ varies  along $\iota$ axis only.
	The common choices for the damping functions are: constant functions, linear functions, quadratic functions, inverse distance functions, etc. \cite{johnson2008notes,kaltenbacher2007numerical}
	
	In order to avoid the long time stability issue in  3D case \cite{kaltenbacher2013modified}, in this example we consider equation (\ref{PML_eq}) in 2D to show the accuracy of the new finite difference scheme.
	In order to obtain the second order temporal accuracy, we introduce the following substitution
	\begin{equation} \label{PML_time_sub}
	u = e^{-\frac{1}{2}\sigma t} \tilde u := s\tilde u.
	\end{equation}
	Then equation (\ref{PML_eq}) can be reformulated as
	\begin{equation}\label{PML_2d}
	\begin{split}
	\frac{1}{\rho c^2} \left(s \tilde u_{tt} + \left(-\frac{1}{4}\sigma^2 + \zeta\right) s \tilde u\right) - \nabla\cdot \left(\frac{1}{\rho}\left(\vec v + \nabla\left(s \tilde u\right)\right)\right) = f, \\ 
	\vec v_t + H\vec v + J \nabla (s\tilde u) = 0. 
	\end{split}
	\end{equation}
	Now we solve equation (\ref{PML_2d}) on the domain $\Omega=[0,2\pi]\times[0,2\pi]$, and $t\in[0,1]$ with the following parameters
	\begin{equation}
	\rho = 1, \ c = 1, \ \sigma_x = \sin(x) - 1, \ \sigma_y = \sin(y) - 1.
	\end{equation}
	Note that  $\rho$ and $c$ are chosen as constants so that the reference solution is available for error calculation. It is worth to mention that the new scheme works well for general cases with variable $\rho$ and $c$. The source term, the initial and boundary conditions are chosen accordingly  to the analytic solution
	\begin{equation}
	u = e^t \sin (x) \sin (y).
	\end{equation}
	
	To demonstrate the convergence order,  the uniform spatial grid size  $h_x = h_y = h_z = h$ and  temporal step size $\tau = \left(\frac{5h}{\pi}\right)^2$ are chosen. The results presented in Table \ref{table_exmaple_4.3} clearly show that the new scheme second-order in time and fourth-order in space, with the truncation error  $O(\tau^2) + O(h^4)$.
	\begin{table}
		\centering
		\caption{Numerical errors for Example (\ref{e.g.:3}) in max norm. Here $\tau = \left(\frac{5h}{\pi}\right)^2$.}
		\label{table_exmaple_4.3}
		\begin{tabular}{|c|c|c|c|c|}
			\hline
			h & $\pi$/25 & $\pi$/50 & $\pi$/75 & $\pi$/100 \\
			\hline
			$\tau$ & 1/25 & 1/100 & 1/225 & 1/400 \\
			\hline
			$E$ & 2.2419e-03 & 1.4182e-04 & 2.8029e-05 & 8.8773e-06 \\
			\hline
			Conv.Order & - & 3.9826 & 3.9986 & 3.9966\\
			\hline
		\end{tabular}
	\end{table}

	\subsection{Example 4}\label{e.g.:4}
	
\begin{figure}[h]
		\caption{(a): acoustic velocity of the Marmousi 2 model. (b): density of the Marmousi 2 model.}
		\vspace{0.5cm}
		\label{ex3_fig1}
		\centering
		\includegraphics[width=1\textwidth]{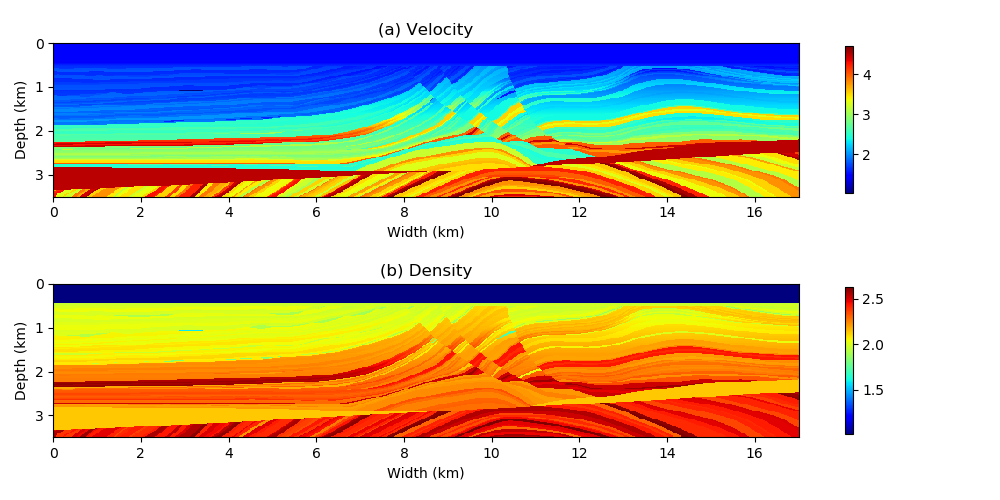}
	\end{figure}

	Finally, we apply the new scheme to  simulate the seismic wave propagation in the Marmousi 2 model \cite{martin2006marmousi2}  using the acoustic wave equation with PML (\ref{PML_eq}).
	The region is a two dimensional domain with depth of $3.5$ km and width of $17$ km, which is discretized using  the spatial step $h_x = h_y = 0.02$ km and time step size $\tau = 0.00125$ s.
	The velocity model and density are  shown in  Figure \ref{ex3_fig1}.
	
\begin{figure}[h]
		\caption{Computation domain with PML zone: $\Omega = \Omega_0 \cup \Omega_{\text{PML}}$, where $\Omega_0$ is the original domain and $\Omega_{\text{PML}}$ is the PML domain with width $0.6$ km}
		\vspace{0.5cm}
		\label{domain_pic}
		\centering
		\includegraphics[width=0.6\textwidth]{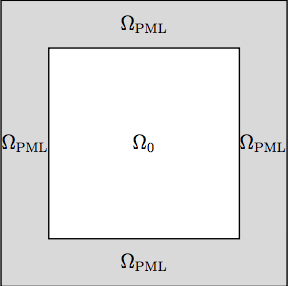}
	\end{figure}

	The Perfect Matched Layer is placed around the model with width $0.6$ km, as depicted in Figure\ref{domain_pic}.
	%
	%
	
	Here the inverse distance damping function is used as:
	\begin{equation}
	\sigma_x(x,z) = \begin{cases}
	0, \quad & \text{in interior}, \\
	\frac{\sigma_{max}}{\text{dist}(x,z)} h_x , \quad & \text{in PML}.
	\end{cases}
	\end{equation}
	Here $\text{dist}(x,z)$ represents the distance of a point $(x,z)$ to the interior domain. The thickness of PML layer is $0.8$ km and $\sigma_{max} = 100$ is used in this example.
	The seismic wave was generated by a Ricker wavelet given by
	\begin{equation}
	s(t,x,z) = \delta(x-x_0,z-z_0) \left(1-2\pi^2 f_p^2 (t-d_r)^2\right) e^{-\pi^2 f_p^2(t-d_r)^2},
	\end{equation}
	where the central frequency $f_p = 5$ Hz and time delay  $d_r=0.2$ s.
	The source  is placed at the centre of the model $x=8.5$ km and $y=1.74$ km to show the absorbing effect of PML.
	
\begin{figure}[h]
		\caption{Snapshots of the wavefield at different time.}
		\label{ex3_fig2}
		\centering
		\includegraphics[width=1\textwidth]{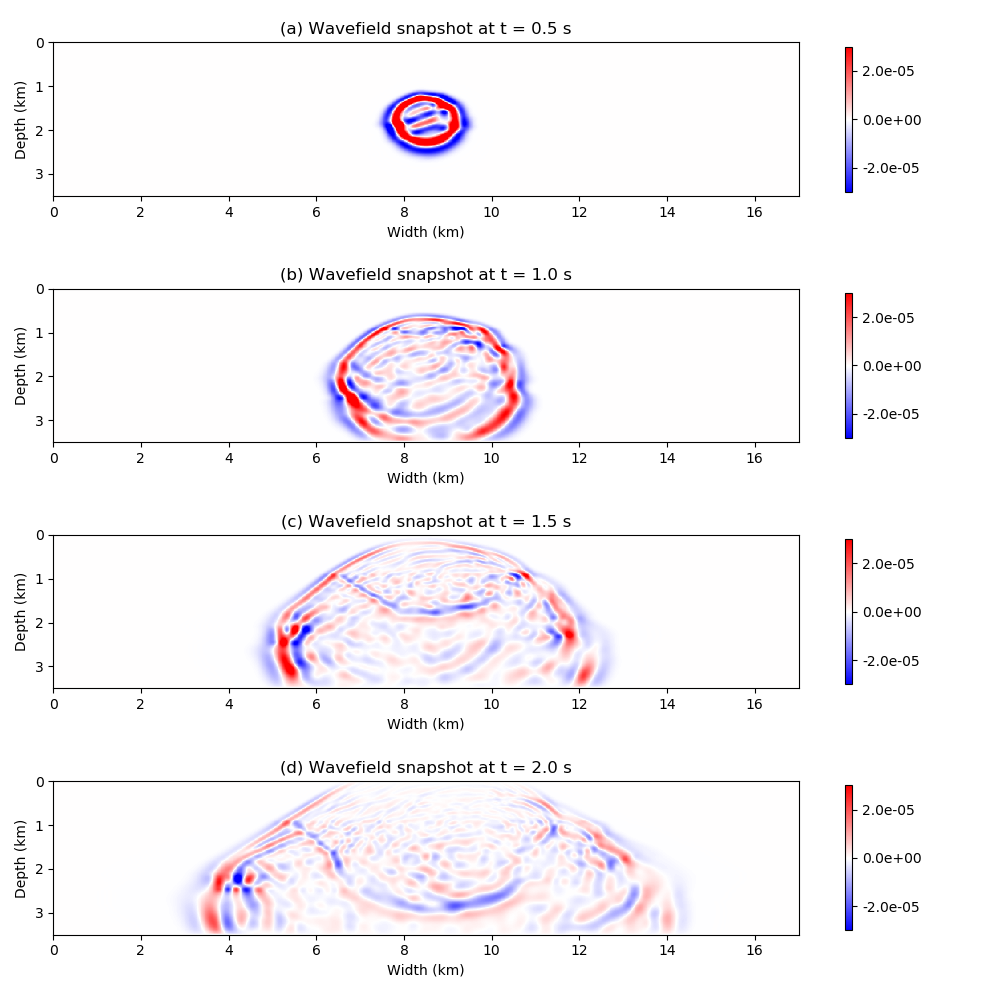}
	\end{figure}

	Four snapshots of the wavefields  are shown in  Figure \ref{ex3_fig2} for  $t=0.5$ s, $1.0$ s, $1.5$ s  and $ 2.0$ s.
	It can be shown in figure (a) that as the seismic wave has not arrived the boundary when $t=0.5$s. Therefore, no wavefields can be observed  except for a small region around the centre of the domain. Then at $t=1.0$s,  one can see in  Figure \ref{ex3_fig2}(b) that the wave approaches the bottom of the domain and the wave energy has been absorbed by the perfectly matched layer.
	Figure \ref{ex3_fig2}(c-d)  show that the seismic wave propagating along the bottom boundary to left and right without  reflection due to PML.
	
\begin{figure}[h]
		\caption{Acoustic energy of the wavefield.}
		\label{ex3_fig3}
		\centering
		\includegraphics[width=0.8\textwidth]{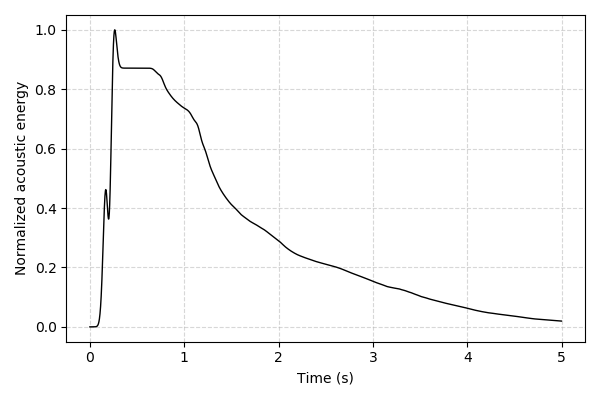}
	\end{figure}

	To furthermore validate the effectiveness of the PML in absorbing energy for reflection reduction,  we calculate the acoustic energy of the acoustic wave as
	\begin{equation}
	E(t) = \int_\Omega \left(\frac{\rho}{2} \vec v \cdot \vec v + \frac{u^2}{2\rho c^2}\right) d x d z.
	\end{equation}
	Here $c$ and $\rho$ are the velocity and density of the Marmousi 2 model, $u$ is the pressure which is computed when  equation (\ref{PML_eq}) is solved. Here  $\vec v$ is the particle velocity which can be computed through the  linear momentum conservation formula
	\begin{equation}\label{energy-acoustic}
	\rho \frac{\partial \vec v}{\partial t} = -\nabla u.
	\end{equation}
	Applying  Crank-Nicolson method to equation(\ref{energy-acoustic}) leads to
	\begin{equation}
	\vec{v}^{n+1} = \vec v^{n} - \frac{\tau}{2\rho}
	\begin{pmatrix}
	\partial_x u^n + \partial_x u^{n+1} \\ \partial_y u^n + \partial_y u^{n+1}
	\end{pmatrix}.
	\end{equation}
	From the energy conservation of wave equation, the acoustic energy should be increasing first and then stay as a constant as time goes on without introducing the PML.
	Figure \ref{ex3_fig3} shows that the acoustic energy inside the domain is decreasing,  which validates that the PML indeed absorbs the energy as seismic wave encounter the boundary of the domain.
	
	\section{Conclusion and future work}
	In this work a new compact explicit finite difference  has been developed to solve acoustic wave equation with variable velocity speed and media density. 
	The second-order convergence in time and fourth-order convergence in space have been theoretically analyzed and numerically confirmed by several examples. The high efficiency of the numerical scheme is obtained through the explicit treatment of the time derivative.  The new scheme has a time complexity which is linear to the total number of grid points for each time step.  Because the spatial differential operator is not self-adjoint, the widely used energy method is not applicable here for stability analysis.  To overcome this issue,  we developed  an empirical 
	stability analysis, which is based on the spectrum estimation of  the discrete differential operator in space, to derive the Courrant-Friedrichs-Lewy condition.  The empirical stability analysis shows that the new scheme is conditionally stable.  Furthermore, we applied the new scheme for the wave equation with PML, which is a more realistic problem in Geophysics. Numerical results from four examples  clearly demonstrated that the new scheme is efficient, stable and accurate for numerical simulation of seismic wave propagation, and  is expected to find wide applications in numerical seismic modelling and related areas.
	
	In the future, the authors plan to generalize this new scheme to more complicated cases such as elastic wave equation. Moreover, some other absorbing  boundary conditions\cite{engquist1977}  will be considered.
	
	\section{Acknowlegements}
	This  work was supported by NSERC (Natural Science and Engineering Research Council of Canada) through  the Discovery
	grant {\bf RGPIN-2019-04830}. The first author also thanks the China Scholarship Council(CSC) for supporting his graduate study at the University of Calgary.

	\appendix
	\section{PML Derivation} \label{appendix_A}
	
	The complex-valued coordinate stretching strategy \cite{chew19943d} has been applied to derive the equation \eqref{PML_eq}.
	Similar derivation can be found in \cite{berenger1994perfectly}.
	To introduce the PML to the 3D acoustic wave equation (\ref{first eq.}), replace the spatial differential operator $\partial_x$, $\partial_y$ and $\partial_z$ by,
	\begin{equation}
	\begin{split}
	\partial_x \longrightarrow \frac{1}{1+\frac{\sigma_x}{\mathrm{i}\omega}}\partial_x = \frac{1}{\eta_x}\partial_x, \\
	\partial_y \longrightarrow \frac{1}{1+\frac{\sigma_y}{\mathrm{i}\omega}}\partial_y = \frac{1}{\eta_y}\partial_y, \\
	\partial_z \longrightarrow \frac{1}{1+\frac{\sigma_z}{\mathrm{i}\omega}}\partial_z = \frac{1}{\eta_z}\partial_z.
	\end{split}
	\end{equation}
	Here $\omega$ is the frequency.
	The damping function $\sigma_x = 0$ in the interior domain and $\sigma_x \neq 0$ in the absorbing layer, and it is invariable along $y$ and $z$ directions.
	The same properties holds for $\sigma_y$ and $\sigma_z$.
	Then the equation (\ref{first eq.}) turns into
	\begin{equation} \label{appendix_1_eq_1}
	\begin{split}
	& \eta_x\eta_y\eta_z \frac{1}{\rho c^2} u_{tt} \\
	-& \left[\left(\partial_x\frac{1}{\rho}\right)\left(\frac{\eta_y\eta_z}{\eta_x}\partial_x u\right) 
	+ \left(\partial_y\frac{1}{\rho}\right)\left(\frac{\eta_x\eta_z}{\eta_y}\partial_y u\right) 
	+ \left(\partial_z\frac{1}{\rho}\right)\left(\frac{\eta_x\eta_y}{\eta_z}\partial_z u\right)\right] \\
	-& \left[\frac{1}{\rho}\partial_x\left(\frac{\eta_y\eta_z}{\eta_x}\partial_x u\right)
	+ \frac{1}{\rho}\partial_y\left(\frac{\eta_x\eta_z}{\eta_y}\partial_y u\right)
	+ \frac{1}{\rho}\partial_z\left(\frac{\eta_x\eta_y}{\eta_z}\partial_z u\right)\right] \\
	& = \eta_x\eta_y\eta_z s
	\end{split}
	\end{equation}
	Here the support of $s$ should be inside the non-absorbing domain, which means $\eta_x\eta_y\eta_z = 1$.
	For the temporal derivative terms in equation (\ref{appendix_1_eq_1})
	\begin{equation} \label{appendix_1_eq_2}
	\begin{split}
	\frac{\eta_x\eta_y\eta_z}{\rho c^2} u_{tt} &= \frac{1}{\rho c^2} \left(1+\frac{\sigma_x}{\mathrm{i}\omega}\right) \left(1+\frac{\sigma_y}{\mathrm{i}\omega}\right) \left(1+\frac{\sigma_z}{\mathrm{i}\omega}\right) u_{tt} \\
	&= \frac{1}{\rho c^2} \left(u_{tt} + \left(\sigma_x+\sigma_y+\sigma_z\right)u_t + \left(\sigma_x\sigma_y+\sigma_x\sigma_z+\sigma_y\sigma_z\right)u + \sigma_x\sigma_y\sigma_x \frac{u}{\mathrm{i}\omega}\right)
	\end{split}
	\end{equation}
	For the spatial derivative terms in equation (\ref{appendix_1_eq_1})
	\begin{equation} \label{appendix_1_eq_3}
	\begin{split}
	\frac{\eta_y\eta_z}{\eta_x}\partial_x u & = \frac{\left(1+\frac{\sigma_y}{\mathrm{i}\omega}\right)\left(1+\frac{\sigma_z}{\mathrm{i}\omega}\right)}{\left(1+\frac{\sigma_x}{\mathrm{i}\omega}\right)} \partial_x u\\
	&= \frac{-\sigma_x + \sigma_y + \sigma_z + \frac{\sigma_y\sigma_z}{\mathrm{i}\omega}}{\mathrm{i}\omega + \sigma_x} \partial_x u + \partial_x u \\
	&:= v_x + \partial_x u.
	\end{split}
	\end{equation}
	For the $v_x$ one has:
	\begin{equation} \label{appendix_1_eq_4}
	\mathrm{i}\omega v_x + \sigma_x v_x + (\sigma_x-\sigma_y-\sigma_z) \partial_x u - \sigma_y\sigma_z \partial\left(\frac{u}{\mathrm{i}\omega}\right) = 0
	\end{equation}
	By using the Fourier transform, one sets $\mathrm{i}\omega v_x = \left(v_x\right)_t$. Then let $w = \frac{u}{i\omega}$, one has $w_t = u$.
	Together with equations (\ref{appendix_1_eq_1}), (\ref{appendix_1_eq_2}), (\ref{appendix_1_eq_3}), (\ref{appendix_1_eq_4}), and replacing the source term $\eta_x\eta_y\eta_z s$ by $s$, the equation system (\ref{PML_eq}) can be derived.

	\nocite{*}
	
	\bibliographystyle{siam}

\begin{thebibliography}{10}

\bibitem{alford1974accuracy}
{\sc R.~Alford, K.~Kelly, and D.~M. Boore}, {\em Accuracy of finite-difference
  modeling of the acoustic wave equation}, Geophysics, 39 (1974), pp.~834--842.

\bibitem{berenger1994perfectly}
{\sc J.-P. Berenger et~al.}, {\em A perfectly matched layer for the absorption
  of electromagnetic waves}, Journal of computational physics, 114 (1994),
  pp.~185--200.

\bibitem{Bilbao2018}
{\sc S.~Bilbao and B.~Hamilton}, {\em Higher-order accurate two-step finite
  difference schemes for the many-dimensional wave equation}, Journal of
  Computational Physics, 367 (2018), pp.~134--165.

\bibitem{britt2018high}
{\sc S.~Britt, E.~Turkel, and S.~Tsynkov}, {\em A high order compact time/space
  finite difference scheme for the wave equation with variable speed of sound},
  Journal of Scientific Computing, 76 (2018), pp.~777--811.

\bibitem{chew19943d}
{\sc W.~C. Chew and W.~H. Weedon}, {\em A 3d perfectly matched medium from
  modified maxwell's equations with stretched coordinates}, Microwave and
  optical technology letters, 7 (1994), pp.~599--604.

\bibitem{chu1998three}
{\sc P.~C. Chu and C.~Fan}, {\em A three-point combined compact difference
  scheme}, Journal of Computational Physics, 140 (1998), pp.~370--399.

\bibitem{cohen1996}
{\sc G.~Cohen and P.~Joly}, {\em Construction analysis of fourth-order finite
  difference schemes for the acoustic wave equation in nonhomogeneous media},
  SIAM Journal on Numerical Analysis, 33(4) (1996), pp.~1266--1302.

\bibitem{engquist1977}
{\sc B.~Engquist and A.~Majda}, {\em Absorbing boundary conditions for
  numerical simulation of waves}, Proceedings of the National Academy of
  Sciences, 74(5) (1977), pp.~1765--1766.

\bibitem{etgen2007computational}
{\sc J.~T. Etgen and M.~J. O’Brien}, {\em Computational methods for
  large-scale 3d acoustic finite-difference modeling: A tutorial}, Geophysics,
  72 (2007), pp.~SM223--SM230.

\bibitem{finkelstein2007finite}
{\sc B.~Finkelstein and R.~Kastner}, {\em Finite difference time domain
  dispersion reduction schemes}, Journal of Computational Physics, 221 (2007),
  pp.~422--438.

\bibitem{horn1991topics}
{\sc R.~A. Horn and C.~R. Johnson}, {\em Topics in matrix analysis, 1991},
  Cambridge University Presss, Cambridge, 37 (1991), p.~39.

\bibitem{johnson2008notes}
{\sc S.~G. Johnson}, {\em Notes on perfectly matched layers (pmls)}, Lecture
  notes, Massachusetts Institute of Technology, Massachusetts, 29 (2008).

\bibitem{kaltenbacher2013modified}
{\sc B.~Kaltenbacher, M.~Kaltenbacher, and I.~Sim}, {\em A modified and stable
  version of a perfectly matched layer technique for the 3-d second order wave
  equation in time domain with an application to aeroacoustics}, Journal of
  computational physics, 235 (2013), pp.~407--422.

\bibitem{kaltenbacher2007numerical}
{\sc M.~Kaltenbacher}, {\em Numerical simulation of mechatronic sensors and
  actuators}, vol.~2, Springer, 2007.

\bibitem{kim1996optimized}
{\sc J.~W. Kim and D.~J. Lee}, {\em Optimized compact finite difference schemes
  with maximum resolution}, AIAA journal, 34 (1996), pp.~887--893.

\bibitem{lele1992compact}
{\sc S.~K. Lele}, {\em Compact finite difference schemes with spectral-like
  resolution}, Journal of computational physics, 103 (1992), pp.~16--42.

\bibitem{levander1988fourth}
{\sc A.~R. Levander}, {\em Fourth-order finite-difference p-sv seismograms},
  Geophysics, 53 (1988), pp.~1425--1436.

\bibitem{li2020efficient}
{\sc K.~Li and W.~Liao}, {\em An efficient and high accuracy finite-difference
  scheme for the acoustic wave equation in 3d heterogeneous media}, Journal of
  Computational Science, 40 (2020), p.~101063.

\bibitem{li2019compact}
{\sc K.~Li, W.~Liao, and Y.~Lin}, {\em A compact high order alternating
  direction implicit method for three-dimensional acoustic wave equation with
  variable coefficient}, Journal of Computational and Applied Mathematics, 361
  (2019), pp.~113--129.

\bibitem{liao2014dispersion}
{\sc W.~Liao}, {\em On the dispersion, stability and accuracy of a compact
  higher-order finite difference scheme for 3d acoustic wave equation}, Journal
  of Computational and Applied Mathematics, 270 (2014), pp.~571--583.

\bibitem{liu2009implicit}
{\sc Y.~Liu and M.~K. Sen}, {\em An implicit staggered-grid finite-difference
  method for seismic modelling}, Geophysical Journal International, 179 (2009),
  pp.~459--474.

\bibitem{martin2006marmousi2}
{\sc G.~S. Martin, R.~Wiley, and K.~J. Marfurt}, {\em Marmousi2: An elastic
  upgrade for marmousi}, The leading edge, 25 (2006), pp.~156--166.

\bibitem{Mattsson2006}
{\sc K.~Mattsson and J.~Nordström}, {\em High order finite difference methods
  for wave propagation in discontinuous media}, Journal of Computational
  Physics, 220(1) (2006), pp.~249--269.

\bibitem{shukla2005derivation}
{\sc R.~K. Shukla and X.~Zhong}, {\em Derivation of high-order compact finite
  difference schemes for non-uniform grid using polynomial interpolation},
  Journal of Computational Physics, 204 (2005), pp.~404--429.

\bibitem{tam1993dispersion}
{\sc C.~K. Tam, J.~C. Webb, et~al.}, {\em Dispersion-relation-preserving finite
  difference schemes for computational acoustics}, Journal of computational
  physics, 107 (1993), pp.~262--281.

\bibitem{wielandt1973eigenvalues}
{\sc H.~Wielandt}, {\em On the eigenvalues of {A+B} and {AB}}, Journal of
  Research of the National Bureau of Standards B, 77 (1973), pp.~61--63.

\bibitem{yang2012central}
{\sc D.~Yang, P.~Tong, and X.~Deng}, {\em A central difference method with low
  numerical dispersion for solving the scalar wave equation}, Geophysical
  Prospecting, 60 (2012), pp.~885--905.

\end{thebibliography}

\end{document}